\documentclass[12pt]{article}
\usepackage{amsmath, amsfonts, amssymb, amsthm} \parskip1pt
\newtheorem{theorem}{Theorem}[section]
\newtheorem{proposition}[theorem]{Proposition}

\newtheorem{corollary}[theorem]{Corollary}
\newtheorem{lemma}[theorem]{Lemma}

\def\T{\mathbb T} 

\def\N{\mathbb N}

\def\frako{{\frak O}}

\def\frako{{\frak O}}

\def\Ga{{\Gamma}}

\def\N{{\bf N}}

\def\ba{{\bf a}}

\def\bB{{\mathbb B}}
\def\bI{{\mathbb I}}

\def\bN{{\mathbb N}}

\def\bQ{{\mathbb Q}}
\def\br{{\bf r}}
\def\bR{{\mathbb  R}}
\def\bv{{\bf v}}

\def\by{{\bf y}}

\def\bz{{\bf z}}
\def\bZ{{\mathbb Z}}

\def\pro{\noindent {\bf{Proof. }}}
\def\card{\rm card}

\date{}

\def\build#1_#2^#3{\mathrel{\mathop{\kern 0pt#1}\limits_{#2}^{#3}}}

\def\smallsquare{\vbox{\hrule\hbox{\vrule height 1 ex\kern 1
ex\vrule}\hrule}}

\def\N{{\mathbb N}}

\def\cE{{\cal E}}
 \def\ocE{\overline {\cE}}
\def\cF{{\cal F}}
\def\cG{{\cal G}}
 \def\ocG{\overline {\cG}}
\def\cS{{\cal E}'}

\def\bp{{\bf p}}
\def\bq{{\bf q}}

\def\bu{{\bf u}}
\def\bv{{\bf v}}
\def\bx{{\bf x}}
\def\by{{\bf y}}
\def\matnm{ {\rm Mat}_{n,m}(\bR)}
\def\bT{{\Theta}}
\def\matnn{ {\rm Mat}_{n,n}(\bR)}

\begin{document}

\title{ Multi-dimensional  metric  approximation by primitive points}
\author{ S. G. DANI, Michel LAURENT,  Arnaldo NOGUEIRA}
\maketitle


\centerline{ABSTRACT} 

\medskip

We consider diophantine inequalities of the form 
$| \bT \bq + \bp - \by |\le \psi(| \bq |)$, with $\Theta \in \matnm$, $\by \in
\bR^n$, where $m,n\in \bN$,  and $\psi$ is a function on $\bN$ with 
positive real values, seeking 
integral solutions  $\bv =(\bq, \bp)^t$ for which the  restriction of $\bv$ 
 to the components of a given partition $\pi$ are primitive integer points. 
In this setting, we establish  metrical statements in the style of the 
Khintchine-Groshev Theorem. Similar solutions
are considered for the doubly metrical inequality $| \bT \bq +\Phi \bp - 
\by |\le \psi(| \bq |)$, with $\Phi \in \matnn$ (other notation as before). 
The results involve the  conditions that $x \mapsto x^{m-1}\psi (x)^n$ 
be non-increasing, and that the components of $\pi$ have at least $n+1$ 
elements each. 

\vskip 7mm

\section{Introduction}
 
For  $d\in \N$ we view $\bR^d$ as the space of $d$-rowed column vectors 
${\bf v}=(v_1,\dots,v_d)^t$, $v_1,\dots,v_d\in \bR$, 
(the $t$ stands for transpose).  
For ${\bf v}\in \bR^d$ denote by $|\bf v |$ the supremum norm in $\bR^d$.  
When ${\bf  w}\in \bR^{m}$ and ${\bf  v}\in \bR^{n}$, we conventionally write as 
$({\bf w,v})^t$
the vector 
$
\left(\begin{matrix}
{\bf  w}
\\
{\bf  v}
\end{matrix}\right)
$
in $\bR^{m+n}$.
 Also, for $m, n\in \bN$, we denote by ${ {\rm Mat}_{n, m}(\bR)}$ the 
vector space of $n\times m$ matrices with real entries. 
  A matrix $X \in {\rm Mat}_{n,m+n}(\bR)$  will be expressed in the form 
$X=(\Theta,\Phi)$ with $\Theta 
 \in { {\rm Mat}_{n,m}(\bR)}$ and $\Phi \in { {\rm Mat}_{n,n}(\bR)}$.
 
 We denote by $P(\bZ^d)$ the set of primitive points in $\bZ^d$, that is, 
the set of   integer $d$-tuples $\bv = (v_1,\dots,v_d)^t$ with $\gcd(v_1,\dots,v_d)= 1$. 
For any subset  $\sigma= \{i_1, \dots , i_\nu\}$ of $\{1,\dots, d\}$ having  $\nu\ge 2$ elements, we denote by 
$P(\sigma)$ the set
 of  integer points $\bv = (v_1,\dots,v_d)^t$ such that $\gcd(v_{i_1},
\dots,v_{i_\nu}) = 1$. Let  $\pi$ be a partition 
of  $\{1, \dots, d\}$ formed by a family of subsets  $\pi_j$  each one with at least two elements.
Then, we denote by $P(\pi)$ the set
 of  integer points $\bv\in \bZ^d$ such that    $\bv \in P(\pi_j)$ for all components  
$\pi_j$. 

Our goal is to refine classical results on metrical diophantine approximation 
for systems of linear or affine inequalities, in the style
of the Khintchine-Groshev Theorem, with  stipulation of additional 
constraints on the solutions being sought, involving 
coprimality; if $\bv$ stands for
the vector of (integral) solution, these conditions are expressed in the form 
$\bv \in P(\pi)$,  $\pi$ a partition of $\{1,\dots , m+n\}$.
The following theorem refines,  in this respect, 
results due to Cassels for $m=1$ and Sprindzuck for $m\ge 2$
 (Theorem II of Chapter VII in \cite{Cas} and Theorem 15 in Chapter I of 
\cite{Sp}); it may be noted however that we need a monotonicity assumption 
 on the approximating function $\psi$, as seen in the statement below, 
though both of the theorems mentioned 
are proved without any monotonicity
condition.

 \begin{theorem}\label{thm1}
 Let $n, m\in \N$  and $\pi$ be a partition of $\{1, \dots, m+n\}$ such 
that every component of $\pi$ has at least $n+1$ elements. 
  Let  $\psi : \bN \to (0,\infty)$ be a  function such that the mapping 
$x \mapsto x^{m-1}\psi(x)^n$ is  non-increasing. If  
 $$
\sum_{j \ge 1} j^{m-1}\psi(j)^n = \infty,
$$
then for almost every pair $(\bT,\by)\in \matnm\times\bR^n$, there exist 
infinitely many points $(\bq,\bp)^t\in P(\pi) $ 
such that
$$
| \bT \bq + \bp - \by |\le \psi(| \bq |). \leqno{(1.1)}
$$
 Conversely, if the series $\sum_{j \ge 1} j^{m-1}\psi(j)^n$ converges, 
then for almost every pair $(\bT,\by)\in \matnm \times\bR^n$, there 
exist only finitely  many integer points $(\bq,\bp)^t\in \bZ^{m+n}$
 for which~$(1.1)$ holds.
 \end{theorem}

In the case $m=n=1$, 
Theorem~\ref{thm1} was proved  in \cite{LaNo},
 under the additional assumption on the growth 
of the function $\psi$, that  $\psi(2 \ell) \gg \psi(\ell)$, namely 
 $\psi(2 \ell)/\psi(\ell)$ bounded below by a positive constant, for all $l\in \bN$; 
see Theorem~2 in \cite{LaNo}.
 
Theorem~\ref{thm1} is a doubly metrical statement in the sense that it concerns 
pairs $(\bT, \by)$. We expect that the same conclusion
 should hold  for any given $\by \in \bR^n$ 
(not just for almost all $\by$) for almost all $\bT$; this is suggested by the
fact that the corresponding statement without the coprimality condition  
$(\bq,\bp)^t\in P(\pi)$  can be deduced,
for every fixed $\by \in \bR^n$, from   the work \cite{ScA} of Schmidt; 
the latter also provides further information, such as counting formulas 
for the number of solutions.  Further metrical results of inhomogeneous approximation
can  also be found in \cite{BuLa}. 
For the case of $\by =0$, viz. the homogeneous case, we are able to establish
such a statement, thus providing a refinement of  the well-known 
Khintchine-Groshev Theorem  \cite{BV,Gro, KlMa}.

\begin{theorem}\label{thm2}
Let $n,m, \pi$ and $\psi$ be as in Theorem~\ref{thm1}. If  
 $
\sum_{j \ge 1} j^{m-1}\psi(j)^n = \infty
$
then for almost every  $\bT \in  \matnm$, there exist infinitely many 
points $(\bq,\bp)^t\in P(\pi) $
such that
$$
| \bT \bq + \bp |\le \psi(| \bq |). \leqno{(1.2)}
$$
 Conversely, if $\sum_{j \ge 1} j^{m-1}\psi(j)^n <\infty$
then for almost every  $\bT\in \matnm $, there exist only finitely  
many $(\bq,\bp)^t\in \bZ^{m+n}$
 for which~$(1.2)$ holds.
 \end{theorem}

Theorems~\ref{thm1} and \ref{thm2} concern  `normalized' systems of
$n$ linear forms  
$\Theta\bq +\bp$ with rank $n$,  meaning that the coefficient matrix
of  $\bp$ is the identity. For general systems with free coefficients, 
we obtain the expected metrical statement  for any fixed $\by$: 

\begin{theorem}\label{thm3}
Let $n,m, \pi$ and $\psi$ be as in Theorem~\ref{thm1}. 
If  $\sum_{j \ge 1} j^{m-1}\psi(j)^n = \infty $
then for every $\by\in\bR^n$ and for almost every  matrix 
$X =(\Theta, \Phi)\in {\rm Mat}_{n,m+n}(\bR) $, there exist infinitely 
many points $(\bq,\bp)^t\in P(\pi) $
such that
$$
| \bT \bq +\Phi \bp -\by |\le \psi(| \bq |). \leqno{(1.3)}
$$
 Conversely, if $\sum_{j \ge 1} j^{m-1}\psi(j)^n<\infty $ 
then for almost every  $X= (\bT,\Phi)\in   {\rm Mat}_{n,m+n}(\bR)$, 
 there exist only finitely many $(\bq,\bp)^t\in \bZ^{m+n}$
 for which $(1.3)$ holds.
 \end{theorem}
 
 It would be worthwhile to illustrate our theorems with two special cases.
 The first one may be viewed as an extension of Theorem 2 of \cite{LaNo} 
to simultaneous approximation, as well as
  an inhomogeneous version of  Gallagher's result \cite{Gal}:
 
 \begin{corollary}
 Let $n$  be a positive integer and let $\psi  : \bN \mapsto \bR^+$ be 
a non-increasing function such that $\sum_{j \ge 1} \psi(j)^n = \infty $.
 For almost every real $2n$-tuple $(\theta_1,\dots ,\theta_n,y_1,\dots , y_n) $ 
there exist infinitely many integer points $(q,p_1,\dots,p_n)$ such that
 $$
 \gcd(q,p_1, \dots, p_n) =1 \quad {\rm and} \quad
 \max_{1\le i \le n}(|q  \theta_i +p_i -y_i |) \le \psi( | q |).
 $$
\end{corollary}

\pro
This is the special case  of Theorem 1.1,  with $m=1$ and $\pi$ the  trivial partition 
of  $\{1,\dots, n+1\}$ (with only one component).
\qed

The following is an application to linear approximation; it may be noted that even 
the homogeneous case with $y=0$ is not obvious when $k\ge 2$. 

 \begin{corollary}
 Let $k$  be a positive integer and let $\psi  : \bN \mapsto \bR^+$ be a 
function such that the 
 mapping $x \mapsto x^{2k-2}\psi(x)$ is  non-increasing and   
$\sum_{j \ge 1} j^{2k-2}\psi(j) = \infty $.
 For every real number $y$ and almost every real linear form $L(\bx) = 
\alpha_1 x_1 +\cdots + \alpha_{2k}x_{2k}$
in $2k$ variables, there exist infinitely many integer points 
$(q_1,\dots,q_{2k})$  such that
 $$
 \gcd(q_{2i-1},q_{2i})=  1 \mbox{ \rm for all }i=1,\dots , k
  \mbox{ \rm and } | L(q_1, \dots, q_{2k}) -y | \le \psi\left( 
\max_{1\le i \le 2k-1}|q_i |  \right) .
  $$
 
\end{corollary}

\pro
This  is the special case of Theorem 1.3, with $m=2k-1, n=1$ and $\pi$ the  
partition as  $\{1,\dots, 2k\}= \coprod_{j=1}^{k}\{2j-1,2j\}$ into pairs of 
indices.
\qed

\medskip

The method of our proofs of Theorems~\ref{thm1} and~\ref{thm2} involves a 
combination of standard 
methods in the metrical theory of diophantine approximation together 
with certain  new  $0-1$ laws. A key observation is that
the set $P(\pi) \subseteq \bZ^{m+n}$ is invariant under  the action  (usual) 
of a certain  subgroup $\Gamma_\pi$ of $SL(m+n,\bZ)$, whose action 
 on ${\rm Mat}_{n,m+n}(\bR)$ by right multiplication is  
ergodic whenever each component of the partition $\pi$ has at least 
$n+1$ elements. 
Corresponding to Theorems 1.1 and 1.2  we define in a natural way certain 
subsets of  ${\rm Mat}_{n,m+n}(\bR)$ which are $\Gamma_\pi$-invariant.
The ergodicity then implies the desired $0-1$ laws for the subsets; the details
involved are described 
in Section~2. The proofs of Theorems 1.1 and 1.2  then proceed along  the  
lines  of that of the Khintchine-Groshev Theorem, 
following the same steps as \cite{BV}.  The main ingredients are adapted 
in Sections~4 and~5  to our present framework involving points in $P(\pi)$. 
 Section~3 is devoted to obtain an  extension of a classical lemma due to 
Cassels that is needed in the proofs of the theorems.  
Theorem~\ref{thm3} readily follows as a consequence of Theorems~\ref{thm1} 
and~\ref{thm2}; it may be worthwhile to note here that alternatively 
Theorem~\ref{thm1} and~\ref{thm2} could be deduced after proving 
Theorem~\ref{thm3} first.

\section{ A general zero-one  law}

Let  $\frak O$ be an infinite subset of $ \bZ^{m+n}$
and  $\psi : \bN \to (0,\infty)$ be a  non-increasing function. We
introduce the following sets depending on  $\frak O$ and $\psi$, for any
$\by \in \bR^n$. For $\bv= (\bq,\bp)^t \in \bZ^{m+n}$ let  us consider the `strip'
$$
S_\bv(\psi,\by) =\Big \{ (\bT,\Phi) \in { {\rm Mat}_{n,m+n}(\bR)} 
: | \bT\bq +\Phi \bp -\by | \le \psi(  | \bq |)\Big\}.
$$ 
For  $(\bT,\Phi) \in { {\rm Mat}_{n,m+n}(\bR)}$ let 
$$\frak O (\bT,\Phi)=\{ \bv \in \frak O : (\bT,\Phi)
\in S_\bv(\psi,\by)\}$$ 
and define  
$$
\cG_\frako  (\psi,\by) =\Big \{ (\bT,\Phi) \in { {\rm
    Mat}_{n,m+n}(\bR)}: \frak O  (\bT,\Phi) \mbox { is infinite}\Big \}. \leqno{(2.1)}
$$
We shall be  concerned with the Lebesgue measure of  
$\cG_{P(\pi)} (\psi,\by)$, and other analogous sets (with $\pi$ as in 
the theorems). 
Towards that end we first consider a related class of sets defined as 
follows. For any  $l \in \N$ let $\psi_l$ be  the function 
defined by $\psi_l (j)=\psi(lj)$ for all $j\in \N$ and set
$$
\cG'_\frako(\psi,\by) =   \bigcap_{l\in \N}\cG_\frako(\psi_l,\by). 
$$
We now describe certain conditions under which the set 
$\cG_\frako' (\psi,\by)$ satisfies the $0-1$ law, namely it is either
 a null set 
(a set of measure $0$)  or a full set (complement  of
a set of  measure $0$), with respect to the Lebesgue measure on 
${\rm Mat}_{n,m+n}(\bR)$. 

In this respect we consider subsets $\frako$ which are orbits of a subgroup
$\Gamma$ of $SL(m+n,\bZ)$. 
We note that on ${\rm Mat}_{n,m+n}(\bR)$ there is an action of 
$SL(m+n,\bZ)$ by (matrix) multiplication on the right. We shall be 
interested in the case when the induced action of $\Gamma$ is ergodic with 
respect to the Lebesgue measure, 
viz. every $\Ga$-invariant 
subset is either a null set or a full set. 

\begin{proposition}\label{prop1}
Let $\frako$ be an orbit of a subgroup $\Gamma$ of $SL(m+n,\bZ)$
whose action on ${\rm Mat}_{n,m+n}(\bR)$  is ergodic, and $\psi$ be
as above. Then 
for every $\by \in \bR^n$,  $\cG'_\frako(\psi,\by) $  
is either a null set or a  full set. 
\end{proposition}

\pro
 We show that $\cG'_\frako(\psi,\by) \bigcap \left(\matnm \times 
GL(n,\bR)\right)$ is invariant under the $\Gamma$-action. Since 
$GL(n,\bR)$  is a full set in $\matnn$ the ergodicity condition 
would then imply that $\cG'_\frako(\psi,\by) $  is either a null 
set or a full set, thus proving  the proposition. 

Let   $(\Theta, \Phi) \in \cG'_\frako( \psi, \by)$, 
with $\Phi \in GL(n,\bR)$, and  $\gamma \in \Gamma$ be given. 
Let  $(\Theta, \Phi)\gamma^{-1}= (\Theta',\Phi')$; to prove the proposition
it suffices to show that 
$(\Theta', \Phi') \in \cG_\frako( \psi_l, \by)$ for all $l\in \N$. 
Let $l\in \N$ be given. 
Choose $a \in \N$  greater than $ (m+n)| \gamma | 
\max(1,2 mn  | \Phi ^{-1} | | \Theta |)$. As  $(\Theta, \Phi)$ belongs 
to $ \cG_\frako (\psi_{a l},\by)$,  there exist  infinitely many  
$\bv=(\bq,\bp)^t  \in \frak O$ such that 
$(\Theta, \Phi) \in S_\bv (\psi_{al}, \by)$.  In other words, we have 
$$
| \bT \bq +\Phi \bp -\by |\le \psi(a l | \bq |)
$$
for these $\bv$. Observe that we have the bounds
$$
| \bp| \le n | \Phi ^{-1}|  | \Phi \bp| \le n | \Phi ^{-1} | 
(| \Theta \bq -\by|
+ \psi (a l | \bq|)) \le n | \Phi ^{-1} | (m | \Theta | | \bq |+| 
\by |+ \psi (a l | \bq|)),$$
and hence $| \bp| \le c| \bq|$
for $c > mn  | \Phi ^{-1} | | \Theta |$, when $| \bq| $ is sufficiently large. 
Now put $\gamma (\bq,\bp)^t = (\bq', \bp')^t$.  We deduce from the above 
upper bounds  the estimate
$$
| \bq'|\le | (\bq', \bp')^t|\le (m+n)| \gamma |  | (\bq, \bp)^t|\le a| \bq|,
$$ 
when $| \bq|$ is sufficiently large. Write 
$$
 \bT'\bq' +\Phi'\bp' =((\Theta, \Phi)\gamma^{-1} )(\gamma (\bq,\bp)^t)
=(\Theta, \Phi) (\bq,\bp)^t = \bT\bq +\Phi\bp.
 $$
  Hence     $ | \bT'\bq' +\Phi'\bp' -\by | =  | \bT\bq +\Phi\bp -\by | 
\le  \psi( a l | \bq |)$. 
As $| \bq' | \leq  a| \bq  | $, and $\psi$ is non-decreasing, we 
have $\psi( a l | \bq |)=\psi_{l}
( a | \bq |)\leq \psi_{l}(  | \bq' |)$, and thus $ | \bT'\bq' +\Phi'
\bp' -\by | \leq  \psi_{l}(  | \bq' |)$. 
Noting that $(\bq',\bp')^t\in \frak O$, this shows that $(\Theta, \Phi) 
\gamma^{-1} =(\Theta', \Phi')\in \cG_\frako(\psi_l,\by)$.
As noted above this proves the  proposition. 
 \qed

\bigskip

We next specialise to the `normalized systems', where  
$\Phi= I_n$, the identity matrix in $\matnn$. Put, with notation $\cG_{\frako}(\psi,\by)$ 
as in~(2.1),
$$
\cE_\frako(\psi,\by) =\Big \{  \Theta \in \matnm  \, : \,  (\Theta,I_n) \in 
\cG_{\frako}(\psi,\by)\Big\}, \leqno{(2.2)}$$
and let 
$$
\cE_\frako(\psi) = \Big\{(\Theta,\by) \in \matnm \times \bR^n : 
\Theta \in \cE_{\frako}(\psi, \by)\Big\};  \leqno{(2.3)} $$
thus $\cE_{P(\pi)}(\psi)$ is the set of all pairs $(\Theta,\by) 
\in \matnm \times \bR^n$ satisfying $(1.1)$. We now define correspondingly
$$
\cS_\frako(\psi,\by) = \bigcap_{l\in \N}\left( \bigcup_{\kappa \in \bN} \cE_\frako(\kappa \psi_l,\by)\right)
\quad {\rm and }\quad
\cS_\frako(\psi) = \bigcap_{l\in \N}\left( \bigcup_{\kappa \in \bN} \cE_\frako(\kappa \psi_l)\right). \leqno{(2.4)}
$$

\begin{proposition}\label{prop2.2}
Let $\frako$ and $\psi$  be  as in  Proposition~\ref{prop1}. Then  $ \cS_\frako(\psi)$  is   
either a null set or a full set.
\end{proposition}

\pro  
Let us introduce the larger set
$$
\cG''_\frako(\psi,\by)  =   \bigcap_{l\in \N}\left( \bigcup_{\kappa \in \bN}
\cG_\frako(\kappa \psi_l,\by)\right)
$$
containing $\cG'_\frako(\psi,\by)$.  Arguing as in the proof of Proposition \ref{prop1},  we observe that 
$\cG''_\frako(\psi,\by) \bigcap \left(\matnm \times  GL(n,\bR)\right)$ is as well invariant under the $\Gamma$-action.
It follows that $\cG''_\frako(\psi,\by) $ is either
a null set or a full set for every fixed $\by \in \bR^n$. Note now that the set of 
$\by\in \bR^n$ for which  $\cG''_\frako(\psi,\by) $ is of full 
Lebesgue measure is invariant under the action of $GL(n,\bR)$ 
by (matrix) left multiplication.  Indeed, if $(\Theta, \Phi) \in
 \cG_\frako(\kappa\psi_l,\by) $,  for any  $g\in GL(n,\bR)$  we have 
 $$
 | g \Theta \bq + g \Phi \bp + g \by | = | g(\Theta \bq + \Phi \bp + \by )| 
\le n | g | | \Theta \bq + \Phi \bp + \by | \le n | g | \kappa \psi_l(| \bq |)
 $$
for infinitely many $(\bq,\bp)^t\in \frako$. Hence $g(\Theta,\Phi)$ belongs to 
$\cG_\frako(\kappa'\psi_l , g\by)$ for  any integer $\kappa' \ge n | g | \kappa$.
 Therefore 
$g\cG''_\frako(\psi,\by)= \cG''_\frako(\psi,g\by)$. 
It follows that  if $\cG''_\frako(\psi,\by) $ is of full Lebesgue measure 
for some $\by \neq 0$ then it is of full Lebesgue measure for all 
$\by \neq 0$. Consequently,  the set 
$$
\cG''_\frako(\psi) :=\{(\Theta, \Phi,\by): \by \in \bR^n,
(\Theta, \Phi) \in \cG''_\frako(\psi,\by) \}
$$ 
is either a full set or a null set  in ${\rm Mat}_{n,m+n}(\bR) 
\times \bR^n.$ Suppose first that it is a full set. Then in particular
 there exists $\Phi \in GL(n,\bR)$ such that
its  fiber  over $\Phi$, namely  
$$
F_\Phi:=\{(\Theta, \by) \in {\rm Mat}_{n,m}(\bR) \times  \bR^n : 
(\Theta,\Phi, \by) \in \cG''_\frako(\psi) \}
$$ 
is a full set in ${\rm Mat}_{n,m}(\bR) \times\bR^n.$ Then 
$\Phi^{-1}F_\Phi=\{(\Phi^{-1}\Theta, \Phi^{-1}\by): (\Theta, \by)\in F_\Phi\}$ 
is also a full set. For $(\Theta_1, \by_1)= (\Phi^{-1}\Theta, \Phi^{-1}\by) 
\in \Phi^{-1}F_\Phi$, and any $(\bq, \bp)^t \in \frak O$ we have 
$$
 | \bT_1\bq + \bp -\by_1 |=  |\Phi^{-1}( \bT\bq + \Phi \bp -\by) 
| \le   n  \vert \Phi^{-1} \vert  | \bT\bq + \Phi \bp -\by|, 
 $$
and since $(\Theta, \Phi) \in \cG''_\frako(\psi,\by) $ this implies that   
$(\Theta_1,\by_1)\in \cS_\frako(\psi)$. We have thus shown that 
$\cS_\frako(\psi)$
contains the full set $\Phi^{-1}F_\Phi$. 
The inverse  argument of multiplication by $\Phi$ show that if 
$\cG''_\frako(\psi)$ is a null set then $ \cS_\frako(\psi) $ is a null 
set. This completes the proof.   \qed

\medskip
An analogous (in fact simpler) argument, with $\cG''_\frako(\psi,{\bf 0})$ 
playing the role of $\cG''_\frako(\psi)$, shows also the following, for 
the set as in $(2.4)$, for $\by ={\bf 0}$; we omit 
the details.  

\begin{proposition}\label{prop2.3}
Let $\Ga$, $\frako$ and  $\psi$   be  as in  Proposition~\ref{prop1}. Then
the set $ \cS_\frako(\psi,{\bf 0})$
 is  either a null set or a full set.
\end{proposition}

\medskip

Before concluding this section  
we introduce the subgroups of $SL(m+n,\bZ)$, and the orbits, to which 
the above results  will be applied in proving the main results. To each partition  
$\pi$ of  $\{1,\dots , m+n\}$ as $ \coprod_{j=1}^k \pi_j $ we associate the
subgroup $\Ga_\pi$ defined as follows. Let $\{e_1,\dots, e_{m+n}\}$ be the
standard basis of $\bR^{m+n}$. For each $j$ let $SL(\pi_j, \bZ)$
be the subgroup of   $SL(m+n,\bZ)$ consisting of elements that fix
$e_i$ for all $i\notin \pi_j$, and let 
 $$
 \Gamma_\pi = \prod_{j=1}^k SL(\pi_j,\bZ). 
 $$
We note that the subgroups $SL(\pi_j, \bZ)$, $j=1,\dots, k$, commute with 
each other
and hence $\Ga_\pi$ is a subgroup of  $SL(m+n,\bZ)$. 
For any $\pi$ as above, we have 
$$
P(\pi) = \Ga_\pi(1,\dots ,1)^t,
$$
an orbit of the $\Ga_\pi$-action on $\bZ^{m+n}$; these
are the orbits $\frako$ involved in the proofs of Theorems~\ref{thm1}, \ref{thm2} and
\ref{thm3}, with the partitions $\pi$ as in the theorems. We will be applying the following 
result on ergodicity, which involves the condition on $\pi$ as in the 
hypotheses of the theorems. 

\begin{proposition}\label{ergod}
   Let $\pi = \coprod_{j=1}^k \pi_j  $ be a partition of $\{1,\dots,
   m+n\}$.  Then the  action of the  
   group $\Gamma_\pi$ on ${\rm Mat}_{n,m+n}(\bR)$, by multiplication
   on the right,  
   is ergodic if and only if for all $j=1,\dots, k$ the  
   cardinality of  $\pi_j$ is at least $ n+1$. 
   \end{proposition}
\pro 
For each $j$ let  $\nu_j$ denote the cardinality of $\pi_j$. Assume first that $\nu_j \ge n+1$ for every $j$ with $1\le j \le k$.
Let $\Omega$ be the open subset of ${\rm Mat}_{n,m+n}(\bR)$  consisting
of all   matrices $M$ such that for each $j$  the $n\times \nu_j$ 
matrix  
$M_{\pi_j}$,  formed by  columns of $M$ corresponding to the indices  in $\pi_j$,
has rank $n$ (the maximum possible).  
Clearly $\Gamma_\pi$ leaves $\Omega$ invariant. Moreover, the
complement of $\Omega$ in  ${\rm Mat}_{n,m+n}(\bR)$ is of Lebesgue
measure $0$ and hence the $\Ga_\pi$-action on ${\rm Mat}_{n,m+n}(\bR)$
is ergodic if and only if the  $\Ga_\pi$-action on $\Omega$ is
ergodic (with respect to the restricted measure). 
Via obvious identifications, we can  write
$
\Omega = \prod_{j=1}^k \Omega_{n, \nu_j}
$
where $\Omega_{n,\nu_j}$  denotes   the open set in  ${\rm
  Mat}_{n,\nu_j}(\bR)$ consisting of matrices of rank $n$. Moreover, 
the action of  
$
 \Gamma_\pi = \prod_{j=1}^k SL(\nu_j,\bZ)
 $
 on $\Omega$ is given by componentwise right multiplication. It follows from Moore's 
  ergodicity theorem (see \cite{Mo}, Theorem~5 and Proposition~6; 
see also \cite{BM}) that the action of  the lattice $SL(\nu_j,\bZ) \subset SL(\nu_j,\bR)$ on 
  $\Omega_{n,\nu_j}$ is ergodic, since the latter can be realised as 
an homogeneous space 
  $  \Omega_{n,\nu_j} =   L_j \backslash SL(\nu_j,\bR)$
  where $L_j$ is the stabiliser of a point in $\Omega_{n,\nu_j}$, and 
$L_j$ is a non-compact closed subgroup.  Thus the action of $ \Gamma_\pi
 $ on the product  $ 
\Omega $ is as well ergodic.  

Suppose now that $\nu_j\leq
n$ for some $j$. Then the determinant of any $\nu_j\times \nu_j$ minor extracted from the component ${\rm Mat}_{n,\nu_j}(\bR)$  remains
invariant under the action of $SL(\pi_j, \bR)$ and hence the action is
not ergodic.
\qed

\section
{Generalisation of a lemma of Cassels}

In this section we prove the following metrical proposition which may be 
viewed as a  multi-dimensional analogue of Lemma~3 in \cite{LaNo}.
The case of $m=1$, which goes back to Cassels \cite{CasB}, deals with boxes  
in $\bR^n$. We are concerned here with a corresponding result for nested 
infinite sequences of product of strips in $\matnm$, viewed as $(\bR^m)^n$.
Such extensions of Cassels' Lemma occur as well in \cite{BVA}.

Let $\Ga$ be a subgroup of $SL(m+n,\bZ)$ and $\frako$ be an infinite 
orbit of $\Ga$ in $\bZ^{m+n}$. For  any  function $\psi$, 
 any $\by\in \bR^n$ and $\Phi \in \matnn$, let 
$$
\cE_{\frako}(\psi,\Phi,\by)  =\Big \{  \Theta \in \matnm  \, : \,  (\Theta,\Phi) 
\in \cG_{\frako}(\psi,\by)\Big\}, \leqno{(3.1)}
$$ 
namely, the  fiber of  $\Phi$ in the set $\cG_\frako(\psi,\by)$ as in~(2.1).

\begin{proposition}\label{cassels}
Suppose that  the function
$\ell\mapsto \ell^{-1} \psi(\ell)$ tends to zero as $\ell$ tends to infinity.
Then for any  $\by\in \bR^n$, $\Phi \in GL(n,\bR)$,   and $\kappa\in \bN$  
the difference 
$$ 
\cE_\frako(\kappa \psi,  \Phi,\by) \setminus \cE_\frako( \psi,\Phi,\by)
$$ 
is  a null set.
\end{proposition}

\pro
For $j=1,\dots, m$ let 
$$
{\frak O}_j=\{(\bq,\bp)^t\in {\frako}: |\bq|=|q_j|\},
$$
  where $q_j$ denotes the $j$th coordinate of $\bq$ and  for any $\Phi \in  GL(n,\bR) $ let  us consider the corresponding set
  $\cE_{\frako_j}(\psi,  \Phi,\by)$.
Clearly 
$$ 
\cE_\frako(\psi, \Phi,\by)=\bigcup_{j=1}^m \cE_{\frako_j}(\psi,  \Phi,\by)
$$
 and hence it suffices to show that  
the difference  $\cE_\frako(\kappa \psi, \Phi,\by) \setminus \cE_\frako( \psi, \Phi,\by)$
 is  a null set for all $j=1,\dots, m$. Moreover, by symmetry
   considerations it suffices to prove this for $j=1$.

We write any $n\times m$ matrix $\Theta \in \matnm$ 
in the form $\Theta = (\xi, \Theta')$, where $\xi$ 
denotes the first column of $\Theta$ and $\Theta' \in {\rm
  Mat}_{n,m-1}(\bR)$ denotes the matrix formed by the remaining columns. 
For any $\bT'\in  {\rm  Mat}_{n,m-1}(\bR)$ let
$$
\cF(\psi, \Theta', \Phi,\by)  = \Big\{\xi \in \bR^n : (\xi, \Theta')
\in \cE_{\frako_1}(\psi, \Phi,\by)\Big\} 
$$
be the fiber over $\Theta'$ of the subset $\cE_{\frako_1}(\psi, \Phi,\by) $.
We write   
$$
\Theta'= (\theta_1, \dots, \theta_n)^t \quad {\rm and}\quad \Phi= (\phi_1, \dots, \phi_n)^t ,
$$
 where $\theta_1, \dots, \theta_n\in \bR^{m-1}$ and $\phi_1, \dots, \phi_n \in \bR^n$ are the 
transposes of the respective rows of
$\bT'$ and $\Phi$. For  any $\bv=(\bq,\bp)^t \in {\frak O}_1$, with  $\bq=(q_1, \bq')^t$ where
$q_1\in \bZ\setminus \{0\}$, $\bq'\in \bZ^{m-1}$,  let
 $\cF_\bv(\psi, \Theta',\Phi,\by)$ be  the set of $\xi \in\bR^n$ for which the matrices $ ((\xi,\Theta'),\Phi)$ belong to $S_\bv(\psi,\by)$. 
It is easily seen that
$$
\cF_\bv(\psi, \Theta',\Phi,\by)  
= \prod_{i=1}^n\left[ {-\phi_i^t\bp- \theta_i^t  \bq'+ y_i
    -\psi(|\bq |)\over | q_1|},{-\phi_i^t\bp- \theta_i^t \bq'+ y_i
    +\psi(|\bq |)\over | q_1|}\right]
$$
is an hypercube in $\bR^n$ and  we can write
$$
\cF(\psi, \Theta',\Phi,\by)  = \limsup_{\bv \in {\frako}_1}\cF_\bv(\psi, \Theta',\Phi,\by) 
$$
as  a limsup of  hypercubes  in $\bR^n$. Observe  that the centers $(-\phi_i^t\bp-\theta_i^t  \bq' + y_i )/ | q_1|$ of
the  intervals occurring in $\cF_\bv(\psi, \Theta',\Phi,\by)$ 
do not depend on $\psi$, and that the length $2 \psi(| \bq |)/ | \bq
|$  is multiplied by the constant 
factor $\kappa$ when $\psi$ is replaced by $\kappa\psi$.  

Let $B$ be a compact subset in $\bR^n$. If for $\bv=(\bq,\bp)^t \in \frako_1$ 
the intersection
$\cF_\bv(\kappa \psi, \Theta',\Phi,\by)\cap B$ is non-empty, then 
we have the rough bound $ | \bp | \le c \kappa | \bq |$ with
$$
c =n | \Phi^{-1} | \left( | \by | + (m-1)  | \Theta'| + \max_{\xi\in B}( | \xi |) + \max_{\ell
  \ge 1}{\psi(\ell)\over \ell}\right),
$$ 
and therefore
$$
\cF(\kappa \psi, \Theta',\Phi, \by) \cap B  = \limsup_{\bv\in 
{\frako_1},\, | \bp | \le c \kappa  | \bq |}
\cF_\bv(\kappa \psi,
\Theta',\Phi,\by)\cap B. 
$$
Now, fix $\kappa\ge 1$ and set
$$
{\frak O}' =\{\bv=(\bq , \bp)^t \in  {\frak O}_1: | \bp | \le c\kappa | \bq | \}. 
$$
Then  we have 
$$
\begin{aligned}
\cF(\psi, \Theta',\Phi,\by) \cap B  & = \limsup_{\bv\in {\frako}'}\cF_\bv(\psi, 
\Theta',\Phi , \by)\cap B, \mbox{ and}
\\
\cF(\kappa \psi, \Theta',\Phi,\by) \cap B  & = \limsup_{\bv\in
  {\frak O}'}\cF_{\bv}(\kappa\psi, \Theta',\Phi, \by)\cap B 
\end{aligned}
$$
If ${\frak O}'$ is finite, both the sets $\cF(\psi, \Theta',\Phi,\by)
\cap B$ and $\cF(\kappa \psi, \Theta',\Phi,\by) \cap B$ are empty. 
If the set ${\frak O}'$ is infinite, we enumerate it as a sequence
$\{\bv_h=(\bq_h, \bp_h)^t\}_{h=1}^\infty$.
Then,  $| \bq_h |$ tends to infinity as $h$ tends to infinity, since
there are 
only finitely many $(\bq, \bp)^t\in {\frak O}'$ with $| \bq |$
bounded. Therefore the length $2 \kappa \psi( | \bq_h|)/ | \bq_h |$ of
the sides of the hypercube 
$\cF_{\bv_h}(\kappa \psi, \Theta',\Phi,\by)$ tends to zero as $h$
tends infinity, since we have assumed that $\lim_{\ell \mapsto
  +\infty} \psi(\ell)/\ell =0$. We may thus  apply the classical lemma 
of Cassels, namely Lemma~9 in \cite{CasB} (see also   Lemma 2.1 in \cite{Har}),  
to the nested sequences of hypercubes 
 $$
 \cF_{\bv_h}(\psi, \Theta',\Phi,\by) \subseteq
 \cF_{\bv_h}(\kappa\psi, \Theta',\Phi,\by), \mbox{ for all } h\ge 1. 
 $$ 
The Lemma asserts that the
associated upper limit sets   
$$
\cF(\psi, \Theta',\Phi, \by) \cap B \,\,  \subseteq \,\,
\cF(\kappa\psi, \Theta',\Phi, \by) \cap B 
$$
have the same Lebesgue measure. Since this holds for any compact set
$B$ it follows  that all the sets $\cF(\kappa\psi,\Theta',\Phi ,
\by), \kappa\ge 1,$ are equal to each other up to a null set. Applying  
Fubini theorem we infer that the fibered sets $\cE_\frako(\kappa\psi, \Phi, \by), 
\kappa\ge 1,$ coincide  up to null sets. \qed

\medskip

The special case $\Phi =I_n$ of Proposition~\ref{cassels} leads to the following 
result for the sets $\cE_\frako (\psi, \by)$ as in $(2.2)$. 

\begin{corollary}\label{cassels2}
Let   $\psi$ be as in Proposition~\ref{cassels}. 
Then for any  $\by\in \bR^n$,  $\kappa\in \bN$ 
the difference  $ \cE_\frako(\kappa \psi,  \by) \setminus 
\cE_\frako( \psi, \by)$ is  a set of  Lebesgue measure~$0$.
\end{corollary}

\section
{Some arithmetical estimates}
In this section we digress and collect some  
 elementary counting results for  use in the estimates needed in the
 next section. 

\begin{lemma}\label{lem:estimate}
Let  $\beta >0 $ and $0<\epsilon <1$ be given. Then for all  
pairs of positive integers $q$ and   $Q$
such that  $q \le Q$, with $Q$ large enough, 
the number of the integers $n$ satisfying
$$
1 \le n \le \beta Q \quad \mbox{and } \quad \gcd(n,q)= 1
$$
 is at least 
$$
 (1-\epsilon)\beta Q \prod_{{  p | q \atop p \,\, { \rm prime} }} \left( 1 - {1\over
     p}\right) =  (1-\epsilon)\beta Q  {\varphi(q)\over q}, 
$$
where $\varphi$ denotes  the Euler totient function.
\end{lemma}

\pro  By the sieving formula of Legendre-Eratosthenes, the number $N$
of these integers $n$ located in the interval $[1,\beta Q]$ and coprime
with $q$ 
is given by the expression
$$
N = \sum_{d | q}\mu(d) \left[ {\beta Q \over d}\right],
$$
where $\mu$ stands as usual for the M\"obius function. Recall that
$\mu(d)$ vanishes unless $d$ is square-free. Therefore 
$$
\begin{aligned}
N  & = \beta  Q \left(\sum_{d | q} {\mu(d)\over   d }\right) - \sum_{d | q}
\mu(d)\left\{{ \beta  Q \over q}\right\} 
\cr
& = \beta  Q \prod_{{  p | q \atop p \,\, { \rm prime} }} \left( 1 -
  {1\over p}\right) \, + \, O\left( 2^{\omega(q)}\right), 
\cr
\end{aligned}
$$
where $\omega(q)$ denotes the number of prime divisors of $q$. Now, it
is well-known that 
$$
\prod_{{  p | q \atop p \,\, { \rm prime} }} \left( 1 - {1\over
    p}\right) ={\varphi(q)\over q} \ge \left( e^{-\gamma} -
  o(1)\right) {1 \over \log\log q }  
\ge   \left( e^{-\gamma} - o(1)\right) {1 \over \log\log Q } ,
$$
where $\gamma$ is the Euler constant, while
$$
\omega(q) \le  \left( 1 + o(1)\right) {\log q \over \log\log q } \le
\left( 1 + o(1)\right) {\log Q \over \log\log Q }. 
$$
Thus, $N$ is asymptotically equivalent to the first term $ \beta  Q
\displaystyle\prod_{{  p | q \atop p \,\, { \rm prime} }} \left( 1 -
  {1\over p}\right) $ as $Q$ tends to infinity, uniformly for  $q$
ranging along  the interval $1 \le q \le Q$.  \qed 

\medskip

\begin{lemma}\label{lem:estimate2}
Let $d$ be an integer $\ge 2$ and $Q$ be a positive real number. Then the 
number of points  
$$
\bq \in P(\bZ^d) \cap \bN^d \quad {\rm with }\quad | \bq | \le Q
$$ 
 is asymptotically equivalent to $\zeta(d)^{-1}Q^d$ as $Q$ tends to infinity.
\end{lemma}

\pro The number of points $\bq$ in $P(\bZ^d)$ with norm $ | \bq | \le
Q$  is asymptotically equivalent to $2^d\zeta(d)^{-1}Q^d$ as $Q$ tends
to infinity; 
this  assertion may be deduced from the special case $K=\bQ$ of
Schanuel's Theorem estimating   
 the number of points in $\mathbb P^{d-1}(K) $ with height bounded by
 $Q$ for any fixed number field $K$ (see for instance  Theorem 5.3 in
 \cite{Lang}). 
 Now, each of the $2^d$ possible quadrants have the same number of
 primitive points contained in the box $| \bq | \le Q$. 
  \qed

Putting  together the two results we next obtain an estimate on certain 
sets of primitive points which will be used in the next section in estimating 
measures of certain sets. 

We first note that given a partition  $\pi = \coprod_{j=1}^{k}\pi_j$ of 
$\{1,\dots ,m+n\}$, where $m,n \in \bN$, by renumbering the indices 
suitably we can arrange so that 
the following holds: there exist indices $a$ and $b$ with $0\le a \le \min(k,m,n)$ and $a \leq b \le k$ such 
that   $\{j, m+j\} \subset \pi_j$ for $1\le j \le a$, 
$\pi_j \subset \{ m+1, \dots,m+n \}$  for $a+1 \le j \le b$, and $\pi_j\subset 
\{ 1, \dots,m \}$ for $b+1\le j \le k$. In other words, the components $\pi_1,\dots , \pi_a$ are those whose intersection with both intervals
$\{1, \dots , m\}$ and $\{m+1, \dots, m+n\}$ is non-empty, while the other components are contained either in $\{1, \dots , m\}$ or in $\{m+1, \dots, m+n\}$.
The case $a=0$ means that there is no component of the first type.

\begin{corollary}\label{cor:estimate}
 Let $m, n \in \bN$  and let $\pi = \coprod_{j=1}^{k}\pi_j$ be 
a partition of $\{1, \dots, m+n\}$ such that each $\pi_j$
has at least $2$ elements. Suppose that the condition as above holds, 
with $0\le a \le b \le k$. For $a+1\leq j \le b$ 
let $d_j$ denote  the cardinality of $\pi_j$. 
Let  $\beta  >0$ and $\epsilon >0$ be given. Then 
for any $\bq =(q_1,\dots, q_m)^t\in \bN^m$ such that $\bq \in P(\pi_j)$ for 
all $j\geq b+1$, and $|\bq| $ is sufficiently large we have 
$$
\card\Big \{ \bp\in \bN^n: \,
(\bq,\bp)^t\in P(\pi), | \bp | \le \beta  | \bq |\Big \}\, 
\, \ge
(1-\epsilon) \beta^{n} | \bq |^n\prod_{j=a+1}^b \zeta (d_j)^{-1}  \prod_{j=1}^a 
{\varphi (q_j)\over q_j}.  
$$
\end{corollary}

\pro For $a+1\le j \le b$  and $\bp \in \bN^n$ 
let $\bp_j$ denote the projection of $\bp$ to the coordinates corresponding
to $\pi_j$ (the latter is contained in $\{m+1, \dots, m+n\}$). Also let $d=\sum_{j=a+1}^b\, d_j$.
Consider the set, say $P'$, 
of $\bp = (p_1,\dots, p_n)^t\in \bN^n$ such that 
(i)~$p_j\leq \beta | \bq | $ and $\gcd (q_j,p_j)=1$ for $1\le j \le a$, (ii)~$  | \bp_j | \le \beta  
| \bq |$ and 
$\bp \in P(\pi_j)$ for $a+1\le j \le b$, and (iii)~$  | \bp | \le \beta  | \bq |$ (the last
part applies afresh only to the remaining $n-a-d$ coordinates of $\bp$ not covered
in (i) and (ii)). 
We note that $P'$ is contained in the set on the left hand side in the 
inequality as above, and the desired cardinality is at least $\card(P')$.  On the other hand  
to pick an element
of $P'$ the choices for $p_1, \dots, p_a$, $\bp_j$, $a+1\leq j \leq b$, 
and the remaining $n-a-d$ coordinates may
be made  independently, satisfying the respective conditions (i),~(ii) and~(iii). 
Therefore when $\bq$ is sufficiently large by Lemmas~\ref{lem:estimate} 
and~\ref{lem:estimate2} we get that $\card(P')$ is bounded below by the 
product of $ (1-\epsilon/3)\beta^a |\bq|^a\prod_{j=1}^a  {\varphi (q_j)\over q_j}$,  
$ (1-\epsilon/3)\beta^{d} |\bq|^{d}\prod_{j=a+1}^b \zeta (d_j)^{-1}$ and 
$ (1-\epsilon/3)\beta^{n-a-d}|\bq|^{n-a-d}$, corresponding to the choices to be made. 
This shows that the estimate as stated in the corollary holds. \qed

\section
{Estimating measures of sets}

\medskip
 Our next objective  will 
be to obtain an estimate as in the following theorem, which
together with the results of the earlier sections will enable us
complete the proofs of Theorems~\ref{thm1}, \ref{thm2} and
\ref{thm3}. Though the components of the partition involved in them are  
assumed to have at least $n+1$ elements, for Theorem~\ref{estimate} 
we need only the  weaker condition as in Corollary~\ref{cor:estimate};
the general form may turn out to be of independent interest. 

Let  
$$
\bI =\left\{\bT\in \matnm : |\bT|\leq \frac 12\right\} 
$$ 
(as before $|\cdot|$ stands 
for the supremum norm).

\begin{theorem}\label{estimate}
 Let $m, n \in \bN$, not both equal to $1$,  and let 
$\pi$ be a partition of $\{1, \dots, m+n\}$ such that each component of 
$\pi$ has at least $2$ elements.
 Let $\lambda$ be the Lebesgue measure on $\matnm$. 
Then there exists a constant $\delta >0$ such that the following
holds: for any  function $\psi : \bN \to (0,\frac 12)$  
such that  the mapping $x \mapsto x^{m-1}\psi(x)^n$ is non-increasing and 
 $$
\sum_{j \ge 1} j^{m-1}\psi(j)^n = \infty, 
$$
and any $\by \in \bR^n$ we have 
$$
\lambda(\cE_{P(\pi)}( \psi,\by ) \cap \bI)\ge \delta.
$$
 \end{theorem}

The proof of this theorem will be completed in the next section. 
In this section we first establish various estimates needed in the proof.  
For simplicity of notation, through the proof we shall suppress $\psi$
from the notation for various sets defined along the way. We note that
the constant $\delta$ is meant to be chosen independently of $\psi$, which 
will be ensured separately in the course of the argument. 

We introduce for any $\by\in \bR^n$ and any $\bv=(\bq,\bp)^t\in
\bZ^{m+n}$  the `strip' 
$$
R_\bv( \by) := \Big \{ \bT \in \matnm; \quad  | \bT\bq  + \bp -\by |
\le \psi(  | \bq |)\Big\}. 
$$
For any  $\bq\in \bZ^m$, let 
$$
\Lambda (\bq)=\{\bv\in
P(\pi) :\bv=(\bq,\bp)^t \mbox { for some }\bp \in \bZ^n\}
\leqno{(5.1)}
$$  
(which could possibly be empty) and  
$$
E_\bq(\by)  = \bigcup_{\bv \in \Lambda (\bq)} R_\bv(\by). 
$$
Since  $\psi( | \bq |) < \frac 12$,   
this is  a  union of disjoint sets.
The fiber $\cE_{P(\pi)}(\psi,\by)$ of $\cE_{P(\pi)}(\psi)$ over $\by \in \bR^n$  is
then equal to  the $ \limsup$ set 
$$
\cE_{P(\pi)}( \psi,\by ) = \bigcap_{Q \ge 1} \bigcup_{ | \bq |  \ge Q} E_\bq(\by).
$$

As usual when dealing with $\limsup$ set in metrical theory, we  first
estimate the Lebesgue   measure of pairwise intersections 
 of the subsets $E_\bq(\by)$, $\bq \in\bZ^m$. 
We  begin by upper bounds. See \cite{BV,Bu,Cas,Gal,Har,ScC,ScA,ScB,Sp,VaB}
for various `overlapping' estimates of this kind 
which are a keystone in metrical diophantine approximation.

\begin{lemma}\label{lem:estimates}
The following estimates hold for any $\by\in \bR^n$: 

{\rm (i)} For any  non-zero  $\bq\in \bZ^m$, 
$$
 \lambda (E_\bq(\by)\cap \bI ) \le 2^n \psi(| \bq |)^n . 
$$

{\rm (ii)} For any  linearly independent vectors $\bq$ and $\bq'$ in
$\bZ^m$,  
$$
\lambda (E_\bq(\by) \cap E_{\bq'}(\by)\cap \bI)  \le 4^n \psi( | \bq | )^n \psi( | \bq' | )^n. 
$$

{\rm (iii)} For $\bq = q\ba$ and $\bq'= q' \ba$, where $q$ and $q'$ are
coprime integers with $| q |\ge | q'| >0$ and  
 $\ba$ is  a non-zero vector in  $ \bZ^m$, we have 
$$
\lambda (E_\bq(\by) \cap E_{\bq'}(\by)\cap \bI)  \le 12^n \psi( | \bq | )^n 
\max\left(\psi( | \bq' | )^n,| q|^{-n}\right). 
$$
\end{lemma}

\pro 
Since we are concerned with upper bounds, it will be convenient to deal  in place of $E_\bq(\by )$ with 
the larger set
$$
    F_\bq(\by )  : =  \bigcup_{\bp\in \bZ^n}
    R_{\bq,\bp)^t}(\by)  
     ={\Big \{} \bT  \in \matnm ;\,  \| \bT\bq  -\by \| \le \psi(  | \bq
  |){\Big \}} ,  
$$
where $\|\cdot\|$ stands, as usual, for the distance to
the nearest point of $\bZ^n$. 
 Let us denote by $\T^n$ the $n$-dimensional torus
$(\bR/\bZ)^n$ and let $\eta :\bR^n \to \T^n$ be the canonical
quotient map. With some obvious abuse of notation, we introduce for any $\by \in \bR^n$, $q\in \bZ\setminus \{0\}$ and $r>0$ the two subsets 
$$
 A_q(\by, r)=\{\bz \in \bR^n : |q\bz-\by|\leq r\} \quad \mbox{ and } \quad 
B_q(\by,r)=\{\bz \in \T^n : \|q\bz-\by\|\leq r\}. 
$$
Thus 
$$
B_q(\by,r) = \bigcup_{\bp \in \bZ^n}\eta( A_q(\by - \bp, r)) = \bigcup_{\bp \in \bZ^n/ q\bZ^n}\eta( A_q(\by - \bp, r)).
$$
noting that  we can restrict the index $\bp$ in the first  union to range along a complete set of representatives modulo $ q $, 
since $ \eta( A_q(\by - \bp, r))=\eta( A_q(\by - \bp_1, r))$ when $\bp \equiv \bp_1$ mod$\,q\bZ^n$.
For any $\bq\in \bZ^m\setminus\{0\}$,  we  introduce the map  $T_\bq:\bI \to \T^n$ 
defined by $T_\bq(\bT)=\eta (\bT \bq)$ for all $\bT\in \bI$, so that 
$$
F_\bq(\by ) \cap \bI = T_\bq^{-1}(B_1(\by, \psi(|\bq |)).
$$

\smallskip
\noindent i) We equip the torus $\T^n$ with the Haar measure $\omega$ normalized by $\omega(\T^n)= 1$ and 
the hypercube $\bI$ with the  Lebesgue measure $\lambda$. 
By a formula of Sprindzuck (see  \cite{Sp}, 
formula $(48)$ on page 35) $T_\bq$ is measure-preserving, and hence 
$$
\lambda(F_\bq(\by)\cap \bI) = \omega(B_1(\by,\psi(| \bq |)) = 2^n \psi(| \bq |)^n,
$$
which proves assertion (i). 

\smallskip
\noindent ii) Let $\by \in \bR^n$ and  $\bq$ and $\bq'$ be linearly
independent vectors in $\bZ^m$. Consider the two maps 
$$
T_\bq : \bI \mapsto  \T^n \quad {\rm and }\quad  T_{\bq'} : \bI
\mapsto \T^n, 
$$
sending $\Theta \in \bI$ to $\eta(\Theta \bq)$ and $\eta(\Theta \bq')$
  respectively. Then, the product mapping 
 $T_\bq\times T_{\bq'} : \bI \mapsto \T^n \times \T^n $ is as well measure-preserving, where 
 $\T^n \times \T^n $ is equipped with the product measure $\omega \times \omega$; this may be seen from 
 formula~$(49)$ in \cite{Sp}. Since 
$$
F_\bq(\by )\cap F_{\bq'}(\by )\cap \bI= (T_\bq\times T_{\bq'})^{-1}\Big(B_1(\by,
\psi(|\bq |))\times B_1(\by,  \psi(|\bq' |))\Big), 
$$
we obtain the equality
$$
\lambda(F_\bq(\by )\cap F_{\bq'}(\by )\cap \bI )=\omega\times\omega \Big(B_1(\by,
\psi(|\bq |))\times B_1(\by,  \psi(|\bq' |))\Big)= 4^n \psi(| \bq
|)^n\psi(| \bq' |)^n, 
$$
which yields (ii).

\smallskip
\noindent iii) Let $\by \in \bR^n$, $\ba\in \bZ^m\setminus\{0\}$, $q,q' \in \bZ$ be two coprime 
integers and let $\bq = q\ba$ and $\bq'= q' \ba$. Then 
$$
F_\bq(\by )\cap F_{\bq'}(\by )\cap \bI= T_\ba^{-1}\Big(B_q(\by, \psi (| \bq |))\cap B_{q'}(\by, 
\psi (| \bq' |))\Big),
$$
Since $T_\ba$ is measure-preserving, we have 
$$
\lambda(F_\bq(\by )\cap F_{\bq'}(\by )\cap \bI) = \omega\Big(B_q(\by, \psi (|\bq |))\cap B_{q'}(\by, 
\psi (|\bq' |))\Big). 
$$

We claim that 
$$
 \omega \Big(B_q(\by, \psi (|\bq |))\cap B_{q'}(\by, 
\psi (|\bq' |))\Big) \le 
 12^n \psi( | \bq | )^n \max\left(\psi( | \bq' | ),{1\over | q|}\right)^n, 
\leqno{(5.2)}
$$
which will  yield (iii).
The set 
$B_q(\by, \psi (|\bq|))\cap B_{q'}(\by, 
\psi (|\bq'|))$ is the image under $\eta$ of its inverse image 
$$ \bigcup_{\bp\in \bZ^n}\bigcup_{\bp'\in \bZ^n} 
\Big(A_q(\by  - \bp, \psi (|\bq|))\cap A_{q'}(\by  - \bp', \psi (|\bq'|))\Big).
$$
For $\bp, \bp'\in \bZ^n$ 
let us set $\alpha (\bp,\bp'): =q\bp'-q' \bp$. The intersection 
$$
A_q(\by - \bp, \psi (|\bq|))\cap A_{q'}(\by - \bp', \psi (|\bq'|))
$$
indexed by the pair $(\bp,\bp')$  is non-empty only 
when the distance between $  q ^{-1}  {(\by  - \bp)}$ 
and $  q'^{-1} {(\by  - \bp') }$
 is at most $ | q |^{-1}{\psi(|\bq |)}+ | q' |^{-1} {\psi(|\bq' |)}$, 
or equivalently when 
$$ 
|\alpha(\bp,\bp') -( q- q')\by|\leq | q' | \psi( | \bq |) + 
| q | \psi( | \bq' |).
$$ 
Observe   that for two integer pairs $(\bp,\bp')$ and $(\bp_1,\bp'_1)$, we have $\alpha (\bp,\bp') =
 \alpha (\bp_1,\bp_1')$   if and only if
 $\bp_1 = \bp +  q \br $ and $ \bp'_1 = \bp' +  q'  \br $ for some $\br \in \bZ^n$,
 since $q$ and $q'$ are coprime integers. Then
 $$
 A_q(\by  - \bp_1, \psi (|\bq|))=A_q(\by - \bp, \psi (|\bq|))  -  \br Ê,\quad  A_{q'}(\by  - \bp'_1, \psi (|\bq'|))=A_{q'}(\by  - \bp', \psi (|\bq'|)) -  \br, 
 $$
so that the associated intersections $A_q(\by - \bp_1, \psi (|\bq|))\cap A_{q'}(\by - \bp_1', \psi (|\bq'|))$ and $A_q(\by - \bp, \psi (|\bq|))\cap A_{q'}(\by -  \bp', \psi (|\bq'|))$
have the same image under $\eta$. Hence, for each $\bz\in \bZ^n$ satisfying
$$
| \bz - ( q - q' )\by | \le | q' | \psi( | \bq |) + | q | \psi( | \bq' |), 
$$
it suffices to keep only one pair $(\bp,\bp')$ with $\alpha(\bp,\bp')= \bz$ included in the union as above.
The number of these integer points $\bz$ is at most $\left( 3 \max( | q' | \psi( | \bq |) + | q | \psi( | \bq' |), 1)\right)^n $
and since $|q'|\le |q|$ it is majorised by 
$$ 
\left( 3 \max( 2 | q | \psi( | \bq' |), 1)\right)^n 
 \le  6^n | q|^n \max\left(\psi( | \bq' | )^n,{| q|^{-n}}\right). 
$$
Also if $\bp, \bp'$ is  such that $\alpha (\bp,\bp')=  \bz $,
then  the Lebesgue measure of the corresponding intersection $A_q(\by  - \bp, \psi (|\bq|))\cap A_{q'}(\by -  \bp', \psi (|\bq'|))$ 
is obviously bounded by
$$
\lambda(A_q(\by  - \bp, \psi (|\bq|)) =  2^n| q |^{-n} { \psi( | \bq |) }^n. 
$$
Together with the preceding 
observation this proves the claim in $(5.2)$ and  hence also (iii). 
\qed

\bigskip

In the opposite direction, we now obtain a lower bound on an average.

\begin{lemma}\label{lem:const}
Let $m,n$ and $\pi$ be as in the statement of Theorem~\ref{estimate}. 
There exists a positive real number $c$, 
such that for any function $\psi$ satisfying the conditions as in 
hypothesis of Theorem~\ref{estimate}, any $\by\in \bR^n$ and any  
sufficiently large integer $Q$, we have the lower bound
  $$
 \sum_{ { \bq \in \bZ^m,\,  1 \le | \bq | \le Q}} \lambda(E_\bq(\by)\cap \bI) \ge  
c \sum_{\ell=1}^Q  \ell^{m-1}\psi(\ell)^n .
  $$

\end{lemma}

\pro
In the sequel, we indicate by $c_1,\dots$ positive real numbers which, 
as well as $c$, depend only on $m$ and $n$ (prima facie they may depend
on $\pi$, but for the latter there are only finitely many possibilities). 

Denote by $| \bq  |_2$ the euclidean norm of the $m$-tuple $\bq$ and 
observe that $|\bq |  \le | \bq |_2$. 
We claim that for all $\bv=(\bq,\bp)^t$ with $\bq\in \bZ^m\setminus\{0\}$, 
$\bp \in \bZ^n$ such that $
| \bp | \le \frac 16 | \bq |_2 
$ we have 
$$
\lambda(R_\bv(\by)\cap \bI ) \ge \left( { \psi( | \bq |)
\over 2^{m-2} (m-1)! | \bq |_2}\right)^n, \leqno{(5.3)}
$$
when $| \bq |$ is sufficiently large. 
 Write the inequality $| \bT\bq + \bp - \by | \le \psi(| \bq |)$ 
in the equivalent form
$$
| \xi_i \cdot \bu- v_i | \le | \bq |_2^{-1}\psi(| \bq | ), \quad 1 
\le  i \le n, \leqno{(5.4)}
$$ 
where $\xi_1, \dots, \xi_n$ are the (column) vectors in $\bR^m$ 
which are the transposes of 
the $n$ rows of $\bT$ and 
$$
\bu = { \bq\over | \bq |_2} , \ \  v_i = {y_i - p_i\over | \bq |_2} 
\quad {\rm for} \quad 1 \le i \le n.
$$
(the dot stands  for the usual scalar 
product in $\bR^m$). Note that $| \bu |_2 =1$. Since $|\bp|\leq \frac 16|\bq|_2$, we have 
$$
\max_{1\le i \le n}| v_i | \le \frac 15, 
$$
 when $| \bq | $ is large enough.
Slicing the euclidean ball $\{\eta \in \bR^m : | \eta |_2 \le \frac 12\}$ 
 by the hyperplanes
 $\eta\cdot\bu = v_i$, $1\leq i \le n$, we see that the set of points 
$\xi_i$ in $\bR^m$
satisfying $(5.4)$ and contained in the hypercube  $| \eta | \le \frac 12$ 
has Lebesgue measure at least
$$
\sigma_{m-1} \left({1\over 4}\right)^{m-1} {2 \psi(| \bq |)\over | \bq |_2} 
\ge { \psi( | \bq |)
\over 2^{m-2} (m-1)! |\bq |_2},
$$
where $\sigma_{m-1}$ denotes the volume of the unit euclidean ball 
in $\bR^{m-1}$. Since $R_\bv(\by)$ is a product set in $(\bR^m)^n$, 
this establishes (5.3).

For any $\bq \in \bZ^m\setminus\{0\}$, denote by $N_\pi(\bq)$ the number of $\bp\in \bZ^n$ 
such that
$$
(\bq,\bp)^t\in P(\pi) \quad \mbox{\rm and} \quad | \bp | \le \frac 16| \bq |. 
$$
It follows from $(5.3)$, with the notation $\Lambda (\bq)$ as before (see 
(5.1)), that 
$$
\lambda(E_{\bq}(\by)\cap \bI) = \sum_{\bv \in \Lambda (\bq)}\lambda(R_\bv(\by) 
\cap \bI) 
\ge N_\pi(\bq)\left( { \psi( | \bq |)\over 2^{m-2} (m-1)! | \bq |_2}\right)^n . \leqno{(5.5)} 
$$

We shall now use the estimate provided by Corollary~\ref{cor:estimate}.
Renumbering the indices if necessary, as in the comment just before the statement of 
the corollary,  we assume that $\pi$ is expressed as  $\pi = \coprod_{j=1}^{k}\pi_j$
so that with  $0\leq a \le b \le k$ the condition formulated there is satisfied.  
Let $$P'(\pi)=\{\bq \in \bN^m : \bq \in P(\pi_j) \mbox { \rm for all } j\geq b+1\}.$$  
Then by Corollary~\ref{cor:estimate} there exists a constant $c_1>0$, 
depending only on $m$ and $n$, such that for any $\bq \in P'(\pi) $ we have 
$$
N_{\pi}(\bq) \ge  c_1 | \bq |^n \prod_{j=1}^a {\varphi (q_j)\over q_j},
$$
and in turn by~(5.5)  
$$ \lambda(E_{\bq}(\by)\cap \bI ) 
\ge  c_2 \psi( | \bq |)^n \prod_{j=1}^{a} {\varphi(q_j)\over q_j},\leqno{(5.6)}
$$
when $|\bq|$ is sufficiently large.

Now  for $\ell\in \bN$ let 
$$
S_\ell= \sum_{{\bq \in P'(\pi), \, | \bq |= \ell}}\,\prod_{j=1}^{a} {\varphi(q_j)\over q_j}.
$$
Then summing the two sides of the  inequality (5.6) over $\bq$, we find
$$
\sum_{\bq \in \bN^m, 1 \le | \bq | \le Q} \lambda(E_\bq(\by)\cap \bI) \ge 
c_2\sum_{\ell=Q_0}^Q S_\ell\psi(\ell)^n \leqno{(5.7)}
$$
for all $Q\ge Q_0$, where $Q_0$ is some  large  integer beyond which 
the above  lower bounds hold.  Assume first that $a\ge 1$. Restricting to $\bq$ with 
$q_1=\ell$ and $1\leq q_i\leq \ell$ for all $i=2,\dots , m$ and using  
Lemma ~\ref{lem:estimate2},  we  obtain the  asymptotical 
lower bound
$$
S_\ell  \ge {\varphi(\ell)\over \ell}\left( \sum_{j=1}^\ell {\varphi(j)
\over j}\right)^{a-1}(c_3\ell^{m-a}),
$$
where  $c_3 = \prod_{j=b+1}^k \zeta(d_j)^{-1}$ with $d_j=\card(\pi_j)$. 
When $a=0$, we combine Lemma~\ref{lem:estimate} for $q=Q=\ell$ with Lemma~\ref{lem:estimate2} to get the  estimate
$
S_\ell  \ge \varphi(\ell) (c_3\ell^{m-2}).
$
Noting that 
the average value of $\varphi(j)/j$ on the interval $[1,\ell]$ is 
asymptotically equal to $1/\zeta(2)$ when $\ell$ is large, we find
that $S_\ell \ge c_4 \varphi(\ell) \ell^{m-2}$ in any case. Hence by  $(5.7)$  
we have
$$
\sum_{ { \bq \in \bZ^m, \,  1 \le | \bq | \le Q}} \lambda(E_\bq(\by)\cap \bI)\ge 
c_5 \sum_{\ell=Q_0}^Q 
{\varphi(\ell)}\ell^{m-2}\psi(\ell)^n. $$
Since the partial sum  $\sum_{j=1}^\ell 
{\varphi(j)\over j}$ is asymptotically equivalent to $\ell/\zeta(2)$ as $\ell$ 
tends to infinity,  and the mapping $\ell\mapsto \ell^{m-1}\psi(\ell)^n$ is 
non-increasing,
by Abel summation we now get that 
$$
\sum_{ { \bq \in \bZ^m, \,  1 \le | \bq | \le Q}} \lambda(E_\bq(\by)\cap \bI)
\ge c_6\sum_{\ell=Q_0}^Q \ell^{m-1}\psi(\ell)^n. 
$$
This completes the proof of the Lemma. \qed

\bigskip
\noindent
{\bf Remark.} In order to bound from below the volume of the intersection of an 
hypercube in $\bR^m$ by a thickened  hyperplane, we have inserted in the hypercube 
an euclidean ball cutting the hyperplane in an $(m-1)$-dimensional ball whose area 
is easily controlled. The recipe  is adequate for our purpose. However, the full 
hyperplane section of the hypercube   is much larger. 
Refined results could eventually be obtained using the explicit formulas of \cite{MaMo, VaB}.

  \section{Completion of  the proof of Theorem~\ref{estimate}}

  We claim that for all $\by \in \bR^n$ we have 
    $$
     \lambda\left( \cE_{P(\pi)}(\psi,\by) \cap \bI \right) \ge
     \delta :=12^{-n}16^{-m}c^2,
    $$
where $c$ is the positive constant as in Lemma~\ref{lem:const}.

Using a   classical  converse to the  Borel-Cantelli Lemma,  we obtain   
the  lower bound 
 $$
 \begin{aligned}
   \lambda\left( \cE_{P(\pi)}(\psi,\by) \cap  \bI\right)&  =  \lambda
\Big( \limsup_{| \bq |\rightarrow + \infty}  
E_\bq(\by)\cap \bI \Big)
   \\
  & \ge \limsup_{Q \rightarrow + \infty}{\left( \sum_{ {1 \le | \bq | \le Q}}
\lambda (E_\bq(\by)\cap \bI)  \right)^2
   \over
  \sum_{ {1 \le  | \bq | \le Q}}\sum_{ { 1 \le  | \bq' | \le Q}} 
 \lambda (E_\bq(\by)\cap E_{\bq'}(\by)\cap \bI ) };
   \end{aligned}
   \leqno{(6.1)}
   $$
    (see for instance  Lemma 2.3 in \cite{Har}).  To find a lower
    bound for the ratio  
as above  we minorize its numerator and majorize the  denominator. For
the numerator,  
Lemma~\ref{lem:const} provides us with the lower bound
    $$
     \sum_{  1 \le  | \bq | \le Q}\lambda ( E_\bq(\by)\cap \bI) 
     \ge c \sum_{ { \ell = 1}}^Q\ell^{m-1}\psi(\ell)^n \leqno{(6.2)}
     $$
 when $Q$ is large enough. Now consider the denominator
    $$
    D : =\sum_{  1 \le  | \bq | \le Q}\,\sum_{  1 \le  | \bq' | \le Q} 
\lambda ( E_\bq(\by)\cap E_{\bq'}(\by)\cap \bI ) .  \leqno{(6.3)}
      $$
When $\bq$ and $\bq'$ are linearly independent, Lemma~\ref{lem:estimates}(ii)
 gives
 $$
 \lambda (E_\bq(\by)\cap E_{\bq'}(\by)\cap \bI )
 \leq 4^n \psi( | \bq | )^n \psi( | \bq' | )^n \le 12^n 
\psi( | \bq | )^n \psi( | \bq' | )^n, \leqno{(6.4)} 
 $$
where the last larger term is introduced for convenience in combining 
with the other terms; see below.  If $\bq, \bq'$ are non-zero linearly 
dependent vectors, say  with 
$| \bq | \ge | \bq' |$,  they can be uniquely written in the form
 $
 \bq = q \ba ,\,\,  \bq' = q' \ba$ with $\ba \in \bZ^m
\setminus\{{\bf 0}\}$, $q\in \bN$, $ q'\in \bZ\setminus\{0\}$, and 
$1 \le | q'
| \le q, \,\, \gcd(q,q')=1.  
 $
Then  by Lemma~\ref{lem:estimates}(iii), 
$\lambda ( E_\bq(\by)\cap E_{\bq'}(\by)\cap \bI )$ is majorised by 
$$ 
12^n \psi( | \bq | )^n \max\left(\psi( | \bq' | ),{1\over  q}\right)^n
\le 12^n \psi( | \bq | )^n \psi( | \bq' | )^n + 12^n 
{\psi( | \bq | )^n\over q^n}. 
$$
 Together with  
(6.3) and (6.4) this yields, on summing the terms corresponding
to all the pairs involved and noting that 
 the number of the integers $q'$ satisfying  $1 \le  | q'| \le q$
 and $ \gcd(q,q')=1$ equals $2 \varphi(q)$, that 
$$
D  \le 12^n\left( \sum_{ { 1\le  | \bq | \le Q}} \psi(| \bq |)^n\right)^2 
+
2\times 12^n \sum_{q=1}^Q{2\varphi(q)\over q^n}\sum_{ {\bq \in q\bZ^m ,1 \le | \bq | \le Q}}\psi(|\bq |)^n. 
  $$ 
Applying the obvious upper  bound $4^m(\ell/q)^{m-1}$ for the 
number of elements $\bv\in\bZ^m$ with norm $| \bv |= \ell/q$,  we
conclude that 
$$
D \le 12^{n}4^{2m}\left(\sum_{ { \ell = 1}}^Q\ell^{m-1}\psi(\ell)^n\right)^2
+12^n 4^{m+1}\sum_{q=1}^Q { \varphi(q)\over q^n}\sum_{ { q\le \ell \le
    Q, \, q | \ell}} 
\left({\ell\over q}\right)^{m-1}\psi(\ell)^n. 
$$
Observe now that
$$
\sum_{ { q\le \ell \le Q, \, q | \ell}}\ell^{m-1}\psi(\ell)^n
\le q^{-1}\sum_{ \ell =1}^{Q}\ell^{m-1}\psi(\ell)^n
$$
since the function $\ell \mapsto \ell^{m-1}\psi(\ell)^n$ is
non-increasing. We thus obtain the bound 
$$
D \le 12^n16^m\left(\sum_{ { \ell = 1}}^Q\ell^{m-1}\psi(\ell)^n\right)^2
+ 12^n 4^{m+1}\left(\sum_{q=1}^Q{ \varphi(q)\over
    q^{m+n}}\right)\left(\sum_{ { \ell =
      1}}^Q\ell^{m-1}\psi(\ell)^n\right).\leqno{(6.5)} 
$$
 Substituting from  $(6.2)$ and $(6.5)$, we get that the right hand
side term in (6.1) is at least 
$$
  \limsup_{ Q \rightarrow + \infty} 
{ \left( c \displaystyle\sum_{ { \ell =
        1}}^Q\ell^{m-1}\psi(\ell)^n\right)^2\over  
  12^n16^m  \left(\displaystyle\sum_{ { \ell =
        1}}^Q\ell^{m-1}\psi(\ell)^n\right)^2+  
12^n 4^{m+1}\left(\displaystyle\sum_{q=1}^Q{ \varphi(q)\over
    q^{m+n}}\right)\left(\displaystyle\sum_{ { \ell =
      1}}^Q\ell^{m-1}\psi(\ell)^n\right) }. $$
We note that since  $m+n\ge 3$,  $\sum_{q=1}^Q{ \varphi(q)\over q^{m+n}}$
 is bounded, while by the condition in the hypothesis $\sum_{ { \ell = 1}}^Q\ell^{m-1}\psi(\ell)^n$ 
can be arbitrarily large when  $Q$ is large enough.  This shows that 
the above limsup is  $12^{-n}16^{-m}c^2$. Therefore  
$ \lambda( \cE_{P(\pi)}(\psi,\by) \cap  \bI )  \ge 12^{-n}16^{-m}c^2,
$
 as claimed. This completes the proof of Theorem~\ref{estimate}. \qed

 \section
 {Proofs of the main results}
 
The fact that the lower bound in Theorem~\ref{estimate}  does not 
depend on the specific function
 $\psi$,  so long as $x\mapsto x^{m-1}\psi (x)^n $ is non-increasing 
and the series 
$\sum_{ { \ell \ge 1}}\ell^{m-1}\psi(\ell)^n$ diverges, is crucial to our 
proofs of Theorems~\ref{thm1} and~\ref{thm2} in the case when $m+n\geq 3$. 
We note that when  the condition holds for a function 
 $\psi$ it also holds  for every scaled version  $\psi_l$, $l \in \bN$, 
(as in \S~2),  since
 $$
 \sum_{ { j \ge 1}}j^{m-1}\psi_l(j)^n= l^{-m+1} \sum_{ { j \ge 1}}(l
 j)^{m-1}\psi(l j)^n 
 \ge l^{-m}\sum_{ { j \ge l}}( j)^{m-1}\psi(j)^n.  
 $$
 
 \subsection
 {Proof of Theorem 1.1}
 Assume first that the series $\sum_{ { \ell \ge 1}}\ell^{m-1}\psi(\ell)^n$ diverges. Put
$$
 \bB= \left\{ \by\in \bR^n : |\by|\leq {1\over 2}\right\}
 $$
and consider the subset $ \cE_{P(\pi)}(\psi_l) \cap (\bI \times \bB)$ of 
$\matnm \times \bR^n$. When $m+n\geq 3$ we get from Theorem~\ref{estimate}
 that there exists $\delta>0$ such that 
 $$
  \lambda\left( \cE_{P(\pi)}(\psi_l,\by) \cap \bI \right) \ge \delta 
  $$
 for all $l \in \bN$ and $\by \in \bR^n$.
 Integrating over 
$\by$ in $\bB$, we obtain the lower bound 
 $$
 \lambda\left( \cE_{P(\pi)}(\psi_l) \cap (\bI \times \bB)\right) 
\ge \delta,\leqno{(7.1)}
$$
for all $l\in \bN$. In the case when $m=n=1$ this statement was 
established in \cite{LaNo} (see bottom of page 422), specifically 
with $\delta =\frac 14$; as noted in the Introduction, 
in \cite{LaNo} Theorem~\ref{thm1} was 
proved for $m=n=1$ under an additional assumption that 
 $\psi(2\ell) \gg \psi(\ell)$; the extra assumption however is not involved
in the proof of the estimate as above. 
Thus (7.1) holds for all $m,n \in \bN$. 
This conclusion in turn implies that  
$$
 \lambda\left( \left( \bigcup_{\kappa \in \bN}\cE_{P(\pi)}(\kappa \psi_l)\right) 
\cap (\bI \times \bB)\right)  \ge \delta.
$$
Hence for  $\cS_{P(\pi)}(\psi)= \bigcap_{l\ge 1}
\left( \bigcup_{\kappa \in \bN}\cE_{P(\pi)}(\kappa \psi_l)\right)$,  
the limit of the decreasing sequence of sets 
$\bigcup_{\kappa \in \bN}\cE_{P(\pi)}(\kappa \psi_l)$, 
we have 
$$
   \lambda\left( \cS_{P(\pi)}(\psi) \cap (\bI \times \bB)\right)  \ge\delta. 
   $$
Thus $\cS_{P(\pi)}(\psi)$ is not a 
null set. Therefore by Proposition~\ref{prop2.2} we get that it  
is a full set. Consequently $\bigcup_{\kappa \in
  \bN}\cE_{P(\pi)}(\kappa \psi)$ is also a full set since it contains $
\cS_{P(\pi)}(\psi)$. Applying   
Corollary~\ref{cassels2} we now obtain finally  that
$\cE_{P(\pi)}(\psi)$ is a full set, proving the first part of
Theorem~\ref{thm1}. 

Assume now that the series $\sum_{\ell \ge 1} \ell^{m-1}\psi(\ell)^n $
converges.  
We have to prove that the set
 $$
\begin{aligned}
 F : =\Big \{ (\bT,\by) \in  \matnm \times  \bR^n : &  
| \bT\bq + \bp -\by |  \le \psi(  | \bq |) 
 \\
& {\rm  for \, infinitely \, many} \,\, (\bq,\bp)^t\in\bZ^{m+n}\Big\}
\end{aligned}
$$
is of Lebesgue measure $0$. The assertion is an easy consequence of 
the Borel-Cantelli Lemma. The case $m=1$ is stated in Theorem~II
in Chapter VII of \cite{Cas}. For completeness, here is a proof  in 
the  general case. Recalling   the sets $F_\bq(\by)$  introduced in  the proof of Lemma 5.2 for any $\bq\in \bZ^m\setminus\{0\}$ and any $\by \in \bR^n$, 
we can write 
$$
F \cap(\bI \times \bR^n) = \coprod_{\by \in \bR^n}\left(\{ \limsup_{\bq\in \bZ^m\setminus\{0\}}F_\bq(\by) \cap \bI \} \times \{\by\}\right),
$$
as a superior limit.  Since $ \lambda(F_\bq(\by)\cap \bI) = 2^n \psi( | \bq |)^n$ for any $\by$, we have 
$$
 \sum_{ { \bq \in \bZ^m, \, 1 \le  | \bq | \le Q}}\lambda (  F_\bq(\by)\cap \bI) 
 = 2^n  \sum_{  \bq \in \bZ^m , \, 1 \le  | \bq | \le Q} \psi( | \bq |)^n
\le 4^m2^n \sum_{ { \ell = 1}}^Q\ell^{m-1}\psi(\ell)^n,
$$
by majorising  the number of points $\bq \in \bZ^m$ with norm
 $| \bq | = \ell$  by $4^m\ell^{m-1}$. 
Now, the Borel-Cantelli lemma yields that $ \displaystyle \limsup_{\bq\in 
\bZ^m\setminus\{0\}}F_\bq(\by) \cap \bI $ is a null set for every $\by \in \bR^n$, 
since the series $\sum_{\ell \ge 1} \ell^{m-1}\psi(\ell)^n $ is convergent. 
Hence $F \cap(\bI \times \bR^n)$ is a null set.
Observe that the subset $F$, naturally embedded in $\bR^{mn+n}$,  is stable 
by the group of integer translations $\bZ^{mn+n}$.
 Therefore, the whole set $F$ is as well a null set. 
This completes the proof of Theorem~\ref{thm1}.

 \subsection
 {Proof of Theorem 1.2}
In the case $m=n=1$ Theorem~\ref{thm2} reduces 
to the classical Khintchine's theorem, see for instance Theorem I in 
Chapter VII of \cite{Cas}; the constraint of coprimality can be readily
met by dividing a (general) solution $(q,p)^t$ by $\gcd(p,q)$, as 
$\psi$ is assumed to be non-increasing. When $m+n\ge 3$ 
the proof of Theorem~\ref{thm2} follows along the same lines as 
the proof of Theorem~\ref{thm1} as above, in fact in a 
much simpler way, with  $\by$ being now 
fixed, equal to the origin ${\bf 0}$ in $\bR^n$; in this case 
Proposition~\ref{prop2.3} 
plays the same role as Proposition~\ref{prop2.2} for Theorem~\ref{thm1}. 
We omit the details.

 \subsection
 {Proof of Theorem 1.3}
 We distinguish the cases $\by \not= 0$ and $\by = 0$. 
  Assume first that $\by $ is a non-zero vector in $\bR^n$. In this case, 
we deduce Theorem~\ref{thm3} from Theorem~\ref{thm1};
  the two statements are in fact equivalent. We have to show that the 
set $\cG_{P(\pi)}(\psi,\by)$ is full (resp. null) 
  for every $\by \in \bR^n\setminus\{0\}$ exactly when the series
   $\sum_{j \ge 1} j^{m-1}\psi(j)^n $ diverges (resp. converges). For this 
purpose, we relate  $\cG_{P(\pi)}(\psi,\by)$ with the various sets
   $$
   \cE_{P(\pi)}(\psi), \quad \cE_{P(\pi)}(\psi,\by), \quad \cE_{P(\pi)}(\psi,
\Phi, \by), \quad \forall  \Phi \in \matnn, \forall\by\in \bR^n,
   $$
   respectively introduced in (2.3), (2.2) and in (3.1). We also set 
   $$
   \cG_{P(\pi)}  (\psi) =\Big \{ ((\bT,\Phi), \by) \in { {\rm
    Mat}_{n,m+n}(\bR)}\times \bR^n:(\Theta,\Phi) \in \cG_{P(\pi)}(\psi,\by)  
\Big \}. 
    $$
 Apart from these sets it would be convenient to introduce certain larger 
sets associated with them,  consisting of the union of the corresponding 
sets over all $\kappa \psi$, as $\kappa$ ranges over $\bN$; 
we shall denote  the corresponding larger set by overlining the notation 
for the original set; thus,  for instance  
   $ \overline {\cE}_{P(\pi)}(\psi) = \cup_{\kappa \in \bN}\cE_{P(\pi)}(\kappa\psi)$. 
This may be compared with the sets introduced in~(2.4). 
   
 \begin{lemma}\label{lem:equivalence}
The following statements are equivalent: 

{\rm (i)} The set $ \ocE_{P(\pi)}(\psi)$ is full (resp. null) in $\matnm 
\times \bR^n$. 

{\rm (ii)} There exists $\Phi \in GL(n,\bR)$ such that  $\ocE_{P(\pi)}
(\psi,\Phi, \by)$ is full (resp. null) in $\matnm$ for almost every 
$\by \in \bR^n$.

{\rm (iii)} For every  $\Phi \in GL(n,\bR)$ the set  $\ocE_{P(\pi)}(\psi,\Phi, \by)$ is full (resp. null) in $\matnm$ for almost every $\by \in \bR^n$.

{\rm (iv)} The set $\ocG_{P(\pi)}(\psi)$ is full (resp. null) in ${\rm Mat}_{n,m+n}(\bR)\times \bR^n$.

{\rm (v)} There exists $\by \in \bR^n\setminus\{0\}$ such that the  set $\ocG_{P(\pi)}(\psi,\by )$ is full (resp. null) in ${\rm Mat}_{n,m+n}(\bR)$.

{\rm (vi)} For every $\by \in \bR^n\setminus\{0\}$  the  set $\ocG_{P(\pi)}(\psi,\by )$ is full (resp. null) in ${\rm Mat}_{n,m+n}(\bR)$.

\end{lemma}

\pro
Arguing as in the proof of Proposition 2.3, observe that the equivariance relations
$$
g\ocE_{P(\pi)}(\psi, \Phi, \by)= \ocE_{P(\pi)}(\psi , g\Phi, g \by) , \quad g\ocG_{P(\pi)}(\psi, \by) = \ocG_{P(\pi)}(\psi, g \by)
$$
hold for any $g\in GL(n,\bR), \Phi \in \matnn , \by \in \bR^n$. Since the action of $GL(n, \bR)$  on $\bR^n\setminus\{0\}$ is transitive, we deduce 
that (v) and (vi) are equivalent. Now  (iv) means that $\ocG_{P(\pi)}(\psi,\by )$ is full (resp. null) for almost all $\by \in \bR^n$ by the Fubini theorem.
Hence (iv) is equivalent to (v) and (vi).
The first equivariance relation  shows in the same way that  (ii) and (iii) are equivalent. Noting that for $\Phi = Id_n$, 
we have $ \cE_{P(\pi)}( \psi,Id_n , \by) = \cE_{P(\pi)}(\psi,\by)$, we obtain the equivalence of (i) with (ii) again by Fubini. 
Finally the equivalence of (iii) and (iv) follows  from the obvious  equality
$$
\ocG_{P(\pi)}( \psi)  = 
\coprod_{\Phi \in \matnn}\coprod_{\by \in \bR^n}\Big(\ocE_{P(\pi)}(\psi, \Phi, \by)\times\{\Phi\}\times \{\by\}\Big),
$$
using Fubini, since $GL(n,\bR)\times \bR^n$ is an open set of  $\matnn\times \bR^n$. \qed
 
Now, Proposition 3.1 tells us that all overlined sets occurring in Lemma 7.1 are equal, up to a null set, to the corresponding
 set without the bar. Thus the same equivalences hold true for the non-overlined sets. In particular,  $ \cG_{P(\pi)}(\psi,\by )$ is full (resp. null)
 if and only if $  \cE_{P(\pi)}(\psi)$ is full (resp. null). By Theorem 1.1, the latter property holds when the series
 $\sum_{j \ge 1} j^{m-1}\psi(j)^n $ diverges (resp. converges).

\medskip

For the case $\by = {\bf 0}$, we use the same strategy, basically multiplying on the left matrices in $ {\rm Mat}_{n,m+n}(\bR)$ or 
in $\matnm$ by $g \in GL(n,\bR)$, to relate the sets 
 $\cG_{P(\pi)}(\psi,{\bf 0})$ and $\cE_{P(\pi)}(\psi,{\bf 0} )$. The proof is simpler and we now apply Theorem 1.2. We omit the details.
  
\medskip

\noindent
{\it Acknowledgements} S.G. Dani and A. Nogueira would like to thank, 
respectively, the Institut de Math\'ematiques de Luminy, Aix-Marseille 
Universit\'e, France and the De Giorgi Center, Italy for hospitality 
while this work was done.

\vskip8mm

SGD : Department of Mathematics, Indian Institute of Technology
Bombay, Powai, Mumbai 400076, India. 

 {\tt sdani@math.iitb.ac.in}

\smallskip

ML \& AN  : Aix Marseille Universit\'e, CNRS, Centrale Marseille, I2M, 
UMR 7373, 13453 Marseille, France 
 
  {\tt michel-julien.laurent@univ-amu.fr, \,\, arnaldo.nogueira@univ-amu.fr}


\begin{thebibliography}{999}

\bibitem{BM} M.B. Bekka and M. Mayer, {\it Ergodic theory and topological 
dynamics of group actions on homogeneous spaces.} London Mathematical 
Society Lecture Note Series, 269. Cambridge University Press, Cambridge, 
2000. x+200 pp. 

\bibitem{BVA}
V. Beresnevich  and  S. Velani,
{\it  A note on zero-one laws in Diophantine approximation}, 
Acta Arithmetica {\bf 133} (2008), 363--374.

\bibitem{BV}
V. Beresnevitch  and  S. Velani,
{\it Classical metric diophantine approximation revisited : the Khintchine-Groshev theorem}, 
Int. Math. Res. Not. IMRN {\bf 2010}, 69--86.

\bibitem{Bu}
Y. Bugeaud,
{\it Approximation by algebraic numbers}, 
Cambridge Tacts in Mathematics, 160. Cambridge University Press, 2004.

\bibitem{BuLa}
Y. Bugeaud and M. Laurent,
{\it Exponents of inhomogeneous Diophantine approximation},
Moscow Math. J. {\bf 5} (2005), 747--766.




\bibitem{Cas}
J. W. S. Cassels,
{\it An Introduction to Diophantine Approximation},
Cambridge Tracts in Math. and Math. Phys., vol. 99, Cambridge
University Press, 1957.

\bibitem{CasB}
J. W. S. Cassels,
{\it Some metrical theorems in Diophantine approximation I},
 Proc. Cambridge Phil. Soc. {\bf 46} (1950), 209--218.



\bibitem{Gal}
P. X. Gallagher,
{\it Metric simultaneous diophantine approximation, II},
Mathematika {\bf 12}   (1965), 123--127.


\bibitem{Gro}
A. Groshev,
{\it Un th\'eor\`eme sur les syst\`emes de formes lin\'eaires},
Dokl. Akad. Nauk SSSR  {\bf 19}   (1938), 151--152.

 \bibitem{Har}
G. Harman,
{\it Metric Number Theory}, 
London Mathematical Society Monographs 18, Oxford University Press, 1998.

\bibitem{KlMa}
D. Y. Kleinbock and G. A. Margulis,
{\it Logarithmic laws for flows on homogeneous spaces},
Invent. Math.  {\bf 138}   (1999), 451--494.



 \bibitem{Lang}
S. Lang,
{\it Fundamentals of Diophantine Geometry}, 
Springer-Verlag, 1983.

\bibitem{LaNo}
M. Laurent and A. Nogueira,
{\it Inhomogeneous approximation with coprime integers and lattice orbits},  
Acta Arithmetica. {\bf 154.4} (2012), 413--427.
 

\bibitem{MaMo}
J.-L. Marichal and M. Mossinghoff,
{\it Slices, Slabs, and Sections of the Unit hypercube},
Online Journal of Analytical Combinatorics, {\bf 3} (2008).




\bibitem{Mo}
C. C. Moore,
{\it Ergodicity of flows on homogeneous spaces}, Amer. J. Math. {\bf
  88} (1966) 154--178.  

\bibitem{ScC}
W. M. Schmidt,
{\it A metrical theorem in diophantine approximation},
Canadian. J.   Math. {\bf 12} (1960) 619--631.  





\bibitem{ScA}
W. M. Schmidt,
{\it Metrical theorem on fractional parts of sequences},
Trans. American.  Math. Soc. {\bf 110} (1964) 493--518.  



\bibitem{ScB}
W. M. Schmidt, 
{\it Diophantine Approximation}, 
Lecture Notes in Math., vol. 785, 1980.

\bibitem{Sp}
V. G. Sprind$\breve{\rm z}$uck, 
{\it Metric Theory of Diophantine Approximations}, 
V. H. Winston \& Sons, 1979.

\bibitem{VaA}
J. D. Vaaler, 
{\it On the metric theory of Diophantine approximation}, Pacific
J. Math. {\bf 76} (1978),   527--539. 


\bibitem{VaB}
J. D. Vaaler, 
{\it A geometric inequality with applications to linear forms},
Pacific J. Math. {\bf 83} (1979),   543--553. 



\end{thebibliography}
\end{document}